\documentclass[12pt,a4paper,reqno,oneside]{amsart}
\usepackage{t1enc}
\usepackage{amssymb}
\usepackage[mathscr]{euscript}
\usepackage{color}
\usepackage{float}
\usepackage{graphicx}
\usepackage{url}
\usepackage{caption}
\usepackage{subcaption}
\usepackage{dsfont}
\usepackage{gensymb}
\usepackage[capitalise]{cleveref}

\numberwithin{equation}{section}

\newtheorem{defi}{Definition}[section]
\newtheorem{thm}[defi]{Theorem}
\newtheorem{lemm}[defi]{Lemma}
\newtheorem{rem}[defi]{Remark}

\newtheorem{cor}[defi]{Corollary}
\newtheorem{prop}[defi]{Proposition}

\newcommand{\TT}{\mathcal{T}}

\newcommand{\M}{\mathcal{M}}

\newcommand{\B}{\mathcal{B}}
\newcommand{\Pp}{\mathcal{P}}
\newcommand{\EE}{\mathcal{E}}

\newcommand{\PP}{\mathbb{P}}

\newcommand{\E}{\mathbb{E}}

\newcommand{\N}{\mathbb{N}}
\newcommand{\R}{\mathbb{R}}
\newcommand{\Z}{\mathbb{Z}}
\newcommand{\Q}{\mathbb{Q}}
\newcommand{\V}{\mathbb{V}}

\newcommand{\GG}{\mathcal{G}}

\newcommand{\FF}{\mathcal{F}}

\newcommand{\1}{\mathds{1}}

\begin{document}

\title[Ergodic Estimation of Dynamic Exceedance Times]{Ergodic Estimation and Model Assessment for Dynamic Exceedance Times}

\author[Å. H. Sande]{Åsmund Hausken Sande}
\address{Å. H. Sande: Department of Mathematics, University of Oslo, Moltke Moes vei 35, P.O. Box 1053 Blindern, 0316 Oslo, Norway.}
\email{aasmunhs@math.uio.no}

\begin{abstract}

{This article concerns the estimation of hitting time statistics for potentially non-stationary processes. The main focus is exceedance times of environmental processes. }

To this end we consider an empirical estimator based on ergodic theory under the assumption that the considered process is a deterministic transformation of some ergodic process. This estimator is empirically analysed and rigorous convergence results, including a central limit theorem, are covered.

Using our estimator, we compute confidence intervals for mean exceedance times of empirical wind data. This serves as a baseline for assessing the performance of several models in terms of predicted mean exceedance time. Special attention is given to the model class known as Gaussian copula processes, which models the environmental process as a deterministic, possibly time-dependent, transformation of a stationary parent Gaussian process.

\end{abstract}

\maketitle

\begin{center}
This Version : \today
\end{center}



\section{Introduction}

The extreme behaviour of environmental processes is an important factor in structural reliability, particularly the study of return periods for extreme events. Consequently there are many different approaches to the modelling of such behaviour in the literature \cite{  jonathan2013statistical, mackay2021effect, papalexiou2018unified, vanemBayes}. 

A popular application of these return periods is through the concept of environmental contours. These are collection of extreme environmental conditions with applications in structural reliability analysis \cite{baarholm2010combining,fontaine2013reliability,giske2018long,vanemTrend}.  For an overview of contour methods see e.g.\ \cite{contcompare,contsummary}. Of particular relevance for this article are convex environmental contours, introduced in \cite{firstaltcontour,altcontour}, and further developed in \cite{dahl2018buffered,voronoi,sande2024convex,sande2023minimal,vanem20193}. These contours are usually constructed based on average hitting times of convex sets by the underlying environmental process. For a few different modelling approaches in this context, we refer to \cite{huseby2023AR1, convcont, mackay2023model, vanem2023analysing, vanem2023analyzing_nonstat, vanemBayes}. Particular recent interest has been given to the incorporation of auto-dependence and nonstationarity in the modelling of the environmental processes. Seasonality is a well know non-stationary effect, but long-term trends are also relevant. For example, long-term trends in significant wave heights have been observed in several analyses \cite{Kushnir,sande2024convex, vanemTrend,vanemBayes}. As such, several approaches to the inclusion of autocorrelation \cite{huseby2023AR1,Leira,mackay2023model,mackay2021effect,vanem2023analysing} and non-stationary behaviour \cite{vanem2023analyzing_nonstat,sande2024convex} have been addressed.

A flexible and convenient class of models are Gaussian copula processes. This class was recently suggested as a universal approach to modelling environmental processes in \cite{papalexiou2018unified}. This class of models have previously been applied to e.g.\ finance \cite{wilson2010copula}, and have been used to estimate return periods of extreme environmental conditions \cite{vanem2023analysing}, including applications to seasonally non-stationary processes \cite{vanem2023analyzing_nonstat}. As such, one of the main goals of this article will be to assess and analyse the performance of this type of model.

To this end we consider an empirical estimator for hitting time statistics for possibly non-stationary processes. The estimator is based on the assumption that our process $V$ is a deterministic, possibly time dependent, transformation of an ergodic process $\widetilde V$. {The main idea of the estimator is the following. Consider the first time $V$ enters some set $\beta$. We can the equivalently consider the first time the ergodic $\widetilde V$ enters a time dependent set $\alpha$. This allows us to extend results from ergodic theory to our original process which provides us with convergence results, including a central limit theorem.}

This estimator can serve as a way to estimate statistical properties of hitting times with a minimal amount of model assumptions. We can also use the empirically estimated return periods as a benchmark to assess the performance of several models in terms of their predicted return periods. Special attention is here given to the aforementioned class of Gaussian copula processes.

{We also analyse various properties of our estimator in order to ensure the validity of our empirical results. This includes the aforementioned convergence results, but also a short study based on synthetic data. The main purpose of this study is to examine our method's inherent bias in the estimation of average exceedance times.}

The article is structured as follows:
{In \Cref{sec:std_ergodics} we give a brief overview of some standard results and definitions needed in the article.} The empirical estimator, and its convergence results, are explained in \Cref{sec:empest}. Results, in their full generality, along with further technical details are given in \Cref{sec:resultsandproofs}. We further analyse the biases and limitations of the empirical estimator in \Cref{sec:statcompare}. Finally, in \Cref{sec:empassessment}, the estimator is applied to empirical wind data, where the predicted return periods of several models is assessed against the baseline of our empirical estimates.













\section{Standard Ergodic Estimators}\label{sec:std_ergodics}

In this {article} we will consider a complete probability space $ \left(\Omega,\FF,\PP\right) $.  We will denote the set of  all non-negative integers, including $ 0 $, by $ \N $, while $\overline\N$ will denote $\N\cup\{\infty\}$.

{The presented estimators will be based on ergodic theory and central limit theory for strong-mixing processes. For an overview we refer to e.g.\ \cite{ergodicbook,mainCLTbook}. We will here mainly consider processes in discrete time.}

\begin{defi}
    By a \textit{discrete process}, we refer to some sequence of random variables $X=\{X_t\}_{t\in\Z}$, where, for some $N\in\N, \, N\geq 1$, we have $X_t(\omega)\in\R^N$ for any $\omega\in\Omega,\, t\in\Z$. 
    We will in this article omit the notation of $\omega$.
\end{defi}

The main result we will use is the mean ergodic theorem. This result holds for any discrete processes that are sufficiently invariant under time shifts.

\begin{defi}
    Define the \textit{shift operator} $S$ as the mapping $ \{X_t\}_{t\in\Z} \mapsto \{X_{t+1}\}_{t\in\Z}$, i.e.\
    $(SX)_t=X_{t+1},$
    for any discrete process $X$. We also denote $(S^sX)_t=X_{t+s},$ for any $s\in\Z$.
\end{defi}

\begin{defi}\label{def:ergodicproc}
    We will refer to a discrete process $X$ as ergodic if the following conditions hold for any $E\in \B(\R^N)^{\otimes\Z}$, the product $\sigma$-algebra of $(\R^N)^{\Z}$.
    $$\PP(X\in E)=\PP(SX\in E),$$
    $$S^sX\in E, \text{ for all $s\in\Z$ } \Rightarrow \PP(X\in E)\in\{0,1\}.$$
\end{defi}

\begin{rem}
    The first condition in \Cref{def:ergodicproc}, is usually referred to as strict stationarity, and can equivalently be stated as
    $$\PP(X_{t_1}\in E_1,
    \cdots, X_{t_n}\in E_n)
    =\PP(X_{t_1+1}\in E_1,
    \cdots, X_{t_n+1}\in E_n),$$
    for all finite sequences, $\{E_i\}_{i=1}^n$, of Borel sets in $\R^N$ and sequences $\{t_i\}_{i=1}^n$, $t_i\in \Z$.
\end{rem}

\begin{thm}[The Mean Ergodic Theorem]\label{thm:mainergodic}
    Let $X$ be an ergodic process and $f$ be any Borel measurable functional $ f:(\R^N)^\Z\to\R $ such that $ \E\left[|f(X) |\right] < \infty$. We then have almost surely that
	\begin{equation*}
		\E\left[f(X) \right]  
		= \lim_{T\to\infty}\frac{1}{T}\sum_{t=0}^{T-1}f(S^tX).
	\end{equation*}
    For a proof, see e.g.\ 
    \cite{ergodicbook}.
\end{thm}

This allows us to compute moments of various functionals. In particular, we will use the fact that hitting times can be considered functionals of the future path of the underlying process.

Additionally, under some additional technical assumptions we can estimate the error of estimates from \Cref{thm:mainergodic} by the following central limit theorem.

\begin{defi}
    For any ergodic process, $X$, we define the \textit{strong mixing coefficient}
    $$\alpha(X,s)=\sup_{A\in\FF_{-\infty}^0,\,B\in \FF_{s}^\infty}\left|\PP(A\cap B)-\PP(A)\PP(B)\right|,$$
    where $\FF_a^b$ is the $\sigma$-algebra generated by $\{X_t\}_{t=a}^b$. The process, $X$ is said to be \textit{strongly mixing} if $\lim_{s\to\infty} \alpha(X,s)=0$.
\end{defi}

\begin{thm}\label{thm:mainCLT}
    Assume $X$ is a strongly mixing process satisfying $\E[X_0^2]<\infty$. We may then define $\M_T=\sum_{t=0}^{T-1}X_t$ and $\sigma_T^2=\V ar[\M_T]$. If $\lim_{T\to\infty}\sigma_T=\infty$ the following are equivalent.
    \begin{itemize}
        \item The family $\{\M_T^2/\sigma_T^2\}_{T\in\N}$ is uniformly integrable.
        \item $(\M_T-\E[\M_T])/\sigma_T$ converges in distribution to a standard normal variable.
    \end{itemize}
    For a proof, see e.g.\ 
    \cite{mainCLTbook}.
\end{thm}

\begin{rem}
    There are many alternative conditions to the uniform integrability of $\{S_T^2/\sigma_T^2\}_{T\in\N}$. A common sufficient criterion ({the Lyapunov condition}) is given in e.g.\ \cite{durrett2019probability}. Here it is assumed that $\E[|X_0|^{2+\delta}]<\infty$ and $\sum_{t=1}^\infty\alpha(X,t)^{\delta/(2+\delta)}<\infty$ for some $\delta>0$. Under this assumption, or similar ones, the central limit theorem will still hold.
\end{rem}

\section{Empirical Estimation of Exceedance Times}\label{sec:empest}

In order to apply ergodic theory to hitting times we will need to slightly modify the estimators to account for censoring. The technical details for the fully general case is given in \Cref{sec:resultsandproofs}, and a simplified version will be given in the remainder of this section.

We now consider a discrete process  $ V:\Omega\times\Z \to \R^N$. In this article $ V $ will model some environmental process measured at all points $ t\in\Z, $ $ 0\leq t < T $ for some final time-point $ T $.

The main quantity of interest is the following.

\begin{defi}
	Assume we have a discrete process $ V $, a time $t\in\Z$, and a set-valued \textit{threshold} function $ \beta: \, \N \to \B(\R^N) $, where $ \B(\R^N) $ denotes the collection of Borel sets of $ \R^N $. We define the \textit{hitting time} of $ V $, with respect to $ \beta $, as the random time, $ \tau : $ $ (\B(\R^N))^\N \times(\R^N)^\Z \times\Z \to \overline\N $ given by
	$$  \tau_\beta(V,t) = \min\{s\geq 0:V_{s+t} \in \beta_s\},  $$
    unless $V_{s+t} \notin \beta_s$ for all $s\geq 0$, in which case we set $\tau_\beta(V,t)=\infty$. Whenever $\beta$ is constant, i.e. $\beta_s=\widetilde\beta$ for all $s\in\N$, we will denote $\tau_\beta(V,t)$ by $\tau_{\widetilde\beta}(V,t)$.
	
\end{defi}
\begin{rem}
	In this article we mainly consider $N=1$ with a constant $ \beta_t=(b,\infty) $. 
    We will, in this setting, refer to $\tau_\beta(V,t)$ as an \textit{exceedance time}.
    
    While we also include arbitrary $N$ for the sake of generality, we will later need a time-varying $\beta$ in order to account for non-stationarity in $V$. 
\end{rem}

\subsection{Stationary Case}\label{subsec:statcase}

{In order to apply our ergodic theory we will in this section assume that the discrete process $\{\tau_\beta(V,t)\}_{t\in\Z}$ is ergodic. Note that ergodicity of $V$ would imply the ergodicity of $\tau_\beta(V,\cdot)$.} This will, in particular, imply that $ \tau_\beta(V,t)=\tau_\beta(V,s) $ in law for any $ t,s\in \Z $. Because of this we may ignore the dependence of $ \tau  $ on $ t $ whenever we are taking expectations. Most importantly, we can apply \Cref{thm:mainergodic} to evaluate $\E\left[f(\tau_\beta(V,0)) \right]$ for any function $ f:\N\to\R $. Mainly, we will be interested in the \textit{return period} $\E\left[\tau_\beta(V,0) \right]$, which can be computed by
\begin{equation}\label{eq:simple_th_estimator}
		\E\left[\tau_\beta(V,0) \right]  
		= \lim_{T\to\infty}\frac{1}{T}\sum_{t=0}^{T-1}\tau_{\beta}(V,t).
\end{equation}

For practical implementations of estimators such as \eqref{eq:simple_th_estimator}, we would want to omit the limit $ T\to\infty $ and use $ T $ equal to the size of our dataset. However, if $\tau_\beta(V,t) \geq T-t$ for any $0\leq t<T$, these estimators would still potentially rely on $V_t$ for $t\geq T$. This constitutes a right-censoring of the data. The easiest way to deal with this is to discard all censored data. For example, we could consider a $\widetilde T \leq T$ such that $\tau_\beta(V,t) < T-t$ for all $t<\widetilde T$ and then use the estimator
\begin{equation}\label{eq:trivialestimator}
    \E\left[\tau_\beta(V,0) \right]  \approx
    \frac{1}{\widetilde T} \sum_{t=0}^{\widetilde T-1}\tau_{\beta}(V,t).   
\end{equation}
This could unfortunately ignore a lot of data, and may consequently yield a higher estimation uncertainty. To remedy this we may replace $\tau_\beta(t)$ with some other quantity which makes use of all the data. Specifically, we will consider

\begin{equation}\label{eq:myestimator}
    \E\left[\tau_\beta(V,0) \right]  \approx
    \frac{1}{T} \sum_{t=0}^{ T-1}\widehat\tau_{\beta}(T;V,t), 
\end{equation}
where
\begin{equation*}
\widehat\tau_\beta(T;V,t) \triangleq 
\begin{cases}%
  \tau_\beta(V,t) & \tau_\beta(V,t) < T-t\\
    T-t+\Lambda(T;V,t) & \tau_\beta(V,t) \geq T-t.
\end{cases}
\end{equation*}
Here $\Lambda:\N\times(\R^N)^\Z\times\N\mapsto\overline\N $ is defined such that $\Lambda(T;V,t)$ depends on $V_t$ only for $0\leq t<T$. Specifically, we will for the rest of the article, except for \Cref{sec:resultsandproofs}, consider the case where
\begin{equation}\label{eq:Lamdadef}
\Lambda(T;V,t) =
\begin{cases}%
    \tau_{S^{T-t}\beta}(V,0) & \tau_{S^{T-t}\beta}(V,0) < T\\
    \infty & \tau_{S^{T-t}\beta}(V,0) \geq T.
\end{cases}
\end{equation}
{Note that we here define $S^s\beta, $ for some $ s\in\N$, by $S^{s}\beta(t)=\beta(t+s)$ for all $t\in\N$.}

\begin{rem}
    {The estimator \eqref{eq:trivialestimator} behaves very similar to \eqref{eq:myestimator}. Consequently, the applications of \Cref{sec:empassessment} could potentially be carried out with either estimator. Nevertheless, \eqref{eq:trivialestimator} appears to perform worse in certain cases. Especially in the case where $\beta_s$ is a decreasing function, i.e.\ $\beta_s\supset \beta_t$ for $s\leq t$. It can be observed from \Cref{prop:nonstatestimator} that this corresponds to the case where $V$ has an increasing trend (see \Cref{subsec:nonstat} for how \eqref{eq:myestimator} is applied to handle a non-stationary $V$). Furthermore, \eqref{eq:myestimator} is monotone in $\beta$ for $\Lambda$ given by \eqref{eq:Lamdadef}. This monotonicity is taken in the sense that $\widehat\tau_\beta(T;V,t) \leq \widehat\tau_\alpha(T;V,t)$ for all $t \leq T$ whenever $\alpha_s \subseteq \beta_s$ for all $s\in\N$.}
\end{rem}

In order to justify our specific choice of $\Lambda$, we consider the following.
Assume that $\tau_\beta(V,t)\geq T-t$ for some $t<T$, this implies
\begin{equation*}
    \tau_{\beta}(V,t)=T-t+\tau_{S^{T-t}\beta}(V,T).
\end{equation*}
Note that if $t\mapsto\tau_\beta(V,t)$ is ergodic then $\tau_{S^{T-t}\beta}(V,T)=\tau_{S^{T-t}\beta}(V,0)$ in law. Our choice of $\Lambda$ then yields that
\begin{equation}\label{eq:conditional_eqlaw}
    \tau_{\beta}(V,t)=\widehat\tau_{\beta}(T;V,t),
\end{equation}
in law for any $0\leq t < T$ conditional on $\tau_{S^{T-t}\beta}(V,0) < T$. Equivalently, this equality in law holds conditional on $\widehat\tau_{\beta}(T;V,t) < \infty$. 


For the sake of implementation, our choice of $ \Lambda $ is equivalent to looping the dataset. We essentially replace $ V_{T+s} $ with $ V_s $ in \eqref{eq:simple_th_estimator}. Specifically, if $\widehat\tau_\beta(V,t)<\infty$ we have
$$\widehat\tau_\beta(T;V,t)=\min\{s\geq 0:V_{t+s \text{ mod } T}\in\beta_s\},
\quad a \text{ mod } b = a-b\lfloor a/b\rfloor.$$

As a consequence of \Cref{exceedfeasiblefull}, in \Cref{sec:resultsandproofs}, we may apply \eqref{eq:myestimator} to achieve a feasible estimator of $\E\left[\tau_\beta(V,T) \right]$. {However, we need a technical requirement on $\Lambda$ to guarantee convergence}. Specifically, we require $\limsup_{T\to\infty}\sup_{0\leq t<T}\Lambda(T;V,t) < \infty$. To resolve this we consider the case where $\beta_s\supseteq\widetilde\beta$ for some measurable $\widetilde\beta\subseteq\R^N$ and for all $s\in\N$, with $\E[\tau_{\widetilde\beta}(V,0)]<\infty$. We then have
\begin{equation*}
    \lim_{T\to\infty}\sup_{0\leq t < T} \Lambda(T;V,t)
    = \sup_{t > 0}\tau_{S^{t}\beta}(V,0)
    \leq \tau_{\widetilde\beta}(V,0)<\infty.
\end{equation*}
almost surely. This additionally guarantees that $\widehat\tau_{\beta}(T;V,t)$ is almost always eventually finite for all $t\leq T$, given a sufficiently high $T$. 

This yields the following results by \Cref{thm:RP_est} (and \Cref{rem:nomomentfinal}). Note that these equalities all hold almost surely.  
\begin{equation}
    \PP\left(\tau_\beta(V,T)> s\right)  
        = \lim_{T\to\infty}\frac{1}{T}\sum_{t=0}^{T-1}\1(\widehat\tau_\beta(T;v,t)> s).\label{eq:simple_cdf}
\end{equation}
If $ \beta_s\supseteq\widetilde{\beta} $ for some measurable $\widetilde\beta\subseteq\R^N$ and all $ s\geq 0 $, with $\E[\tau_{\widetilde\beta}(V,T)^2]<\infty$ then
\begin{equation}
    \E\left[\tau_\beta(V,T)\right] 
        = \lim_{T\to\infty}\frac{1}{T}\sum_{t=0}^{T-1}\widehat\tau_\beta(T;V,t),\label{eq:simple_mean}
\end{equation}
If we further assume that $\E[\tau_{\widetilde\beta}(V,T)^3]<\infty$ then
\begin{equation}
    \E\left[\tau_\beta(V,0)\tau_\beta(V,i) \right] 
        = \lim_{T\to\infty}\frac{1}{T}\sum_{t=i}^{T-1}\widehat \tau_\beta(T;V,t-i)\widehat\tau_\beta(T;V,t),\label{eq:simple_cov}
\end{equation}

In the case where $ \beta $ is constant, i.e.\ $ \beta_s = \widetilde{\beta} $ for all $ s $, we can derive a computationally cheap explicit formula for \eqref{eq:simple_mean}. First we define $ \Theta $, the set of all observed points where the process is in $ \widetilde{\beta} $, i.e. $ \Theta\triangleq\{\theta\in\N: \theta<T,\, V_\theta \in\widetilde{\beta}\} $. We assume that the size of this set, $ |\Theta| $, is at least $ 1 $. We denote by $ \theta_i $, $ i\leq|\Theta| $ the $ i $'th element of $ \Theta $, and define $ \theta_{|\Theta|+1}=T+\theta_1 $. We lastly introduce the shorthand $\Delta\theta_i=\theta_{i+1}-\theta_i$. This yields
\begin{equation}\label{eq:statmean}
  \frac{1}{T}\sum_{t=0}^{T-1} \widehat\tau_\beta(t)
	=\sum_{i=1}^{|\Theta|}
            \frac{\Delta\theta_i(\Delta\theta_i-1)}{2 T}.  
\end{equation}

{Finally, in \Cref{thm:feasibleCLT}, we also prove a central limit theorem for our estimator. Specifically,}
\begin{equation}\label{eq:simpleCLT}
    \frac{ \frac{1}{T}\sum_{t=0}^{T-1} \widehat\tau_\beta(V,t) - \E[\tau_\beta(V,0)]}{{\text{sd}\left(\frac{1}{T}\sum_{t=0}^{T-1} {\tau_\beta}(V,t)\right)}},
\end{equation}
converges in distribution to a standard normal random variable when $T\to\infty$. In order to practically compute this we note that
\begin{equation*}
    \text{sd}\left(\sum_{t=0}^{T-1} {\tau_\beta}(V,t)\right)^2 
    =
    T \,\E[\tau_\beta(V,0)^2]+
    2\sum_{i=1}^{T-1}\E\big[\tau_\beta(V,0)\tau_\beta(V,i)\big](T-i),
\end{equation*}
which is estimable from \eqref{eq:simple_cov}.

\subsection{Non-Stationary Case}\label{subsec:nonstat}
{

The methods presented so far require a stationary $ V $. This usually ignores several important factors, mainly seasonal behaviour and long-scale trends. In order to include and examine such effects we will need to slightly adjust our empirical estimator. 

\begin{prop}\label{prop:nonstatestimator}
    Consider a discrete, possibly non-stationary, process $V$ and a threshold function $ \beta: \, \N \to \B(\R^N) $. Let $ G:\Z \times \R^N \to \R^M $, for some $M\in\N$, be any measurable function such that $v \mapsto G_t(v)$ is injective for any $t\in\Z$. 
    
    We define $\widetilde{V}$ by $ \widetilde{V}_t= G_t(V_t)$ for all $t\in\Z$. For any $ u\in\Z $ we define $ \alpha_s=G_{u+s}(\beta_s) $. We then have
	
	$$\tau_{{\beta}}(V,u)=\tau_{\alpha}(\widetilde{V},u).$$
	
	\begin{proof}
		Follows from the observation that	
		$$ V_{u+s}\in{\beta_s} \Leftrightarrow 
		G_{u+s}(V_{u+s})\in G_{u+s}({\beta})  \Leftrightarrow 
		\widetilde{V}_{u+s}\in\alpha_s.$$
	\end{proof}
\end{prop}

If $\widetilde V$ in \Cref{prop:nonstatestimator} is ergodic we can apply the theory from \Cref{subsec:statcase}.
This means we can estimate the average hitting time of $ V $ with respect to $\beta$ by applying our estimator to $ \widetilde{V} $ and $\alpha$.   Specifically, we have almost surely that
\begin{equation*}
    \E\left[\tau_\beta(V,u) \right]  =
    \lim_{T\to\infty} \frac{1}{T} \sum_{t=0}^{T-1}\widehat\tau_{\alpha}(T;\widetilde V,t),
\end{equation*}
if the technical conditions of \eqref{eq:simple_mean} hold. In this context we would require $\alpha_s\supseteq\widetilde{\alpha} $ for some measurable $\widetilde\alpha\subseteq\R^M$ and all $ s\geq 0 $, with $\E[\tau_{\widetilde\alpha}(\widetilde V,T)^2]<\infty$.
}

\section{Analysis of Bias for Synthetic Data}\label{sec:statcompare}

In this section we aim to examine the bias and limitations of the estimator by means of synthetic data. This will also allow us to comment on the benefits of including non-stationarity in the modelling of $V$.

For the sake of clarity and simplicity, we will consider the average first exceedance time of a one-dimensional and ergodic $V$ above a threshold $b$. Specifically, we will consider the functions
\begin{equation*}
    \EE_b(V) \triangleq \E[\tau_{(b,\infty)}(V,0)], \quad 
    \widehat\EE_b(V,T) \triangleq \frac{1}{T}\sum_{t=0}^{T-1}\widehat\tau_{(b,\infty)}(V,t).
\end{equation*}
Note that the random variable $\widehat\EE_b(V,T)$ is taken to be infinite if $\widehat\tau_{(b,\infty)}(V,t)=\infty$ for any $t$.

We may note from \eqref{eq:conditional_eqlaw} that the estimator is unbiased conditional on the event event that it is finite, i.e.\ 
$$\E\left[\tau_{(b,\infty)}(V,0) -\widehat\EE_b(V,T) \Big| \widehat\EE_b(V,T)<\infty \right]=0.$$ 
Additionally, we know that for a fixed $b$ we have $\lim_{T\to \infty}\widehat\EE_b(V,T) = \EE_b(V)$.
However, for a dataset of fixed size $T$, if we compute the curve $b\mapsto\widehat\EE_b(V,T)$ for all $b$ such that $\widehat\EE_b(V,T)<\infty$,  we introduce a bias. This bias comes from the fact that 
$$\EE_b(V) -\E\left[\widehat\EE_b(V,T) \Big| \widehat\EE_b(V,T)<\infty \right] \geq 0,$$
where the inequality is usually strict. Fortunately, as we will see, this bias vanishes for smaller values of $b$.

\subsection{Stationary Case}

To examine this effect we consider $V$ as a standard normal AR(1) process. We can therefore define $V$ by 
\begin{equation}\label{eq:AR1_V}
    V_t=\sqrt{1-\rho^2}\sum_{s=-\infty}^t \rho^{t-s}X_s=\rho V_{t-1} + \sqrt{1-\rho^2}X_t,
\end{equation}
for all $t\in\Z$. Here $\rho=97\%$ and $\{X_t\}_{t\in\Z}$ is an independent sequence of standard normal random variables. Note that we could here equivalently define $V$ as the unique standardised (mean 0 and variance 1) Gaussian process with autocorrelation $\rho(V,t)=\text{cor}(V_0,V_t)=\rho^t$.
\begin{rem}
    Note that since $\EE_b(sV+m)=\EE_{(b-m)/s}(V)$, there is no loss of generality from considering a standardised process.
\end{rem}
\begin{figure}[h]
    \centering
    \includegraphics[width=0.75\textwidth]{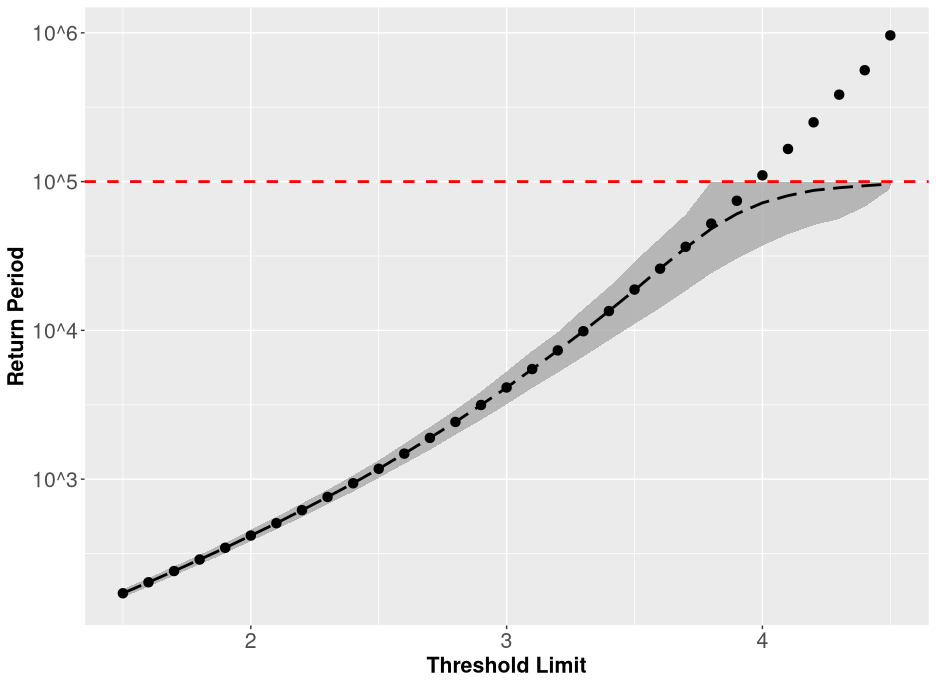}
    \caption{Plot of conditional mean (black dashed line) and $80\%$ prediction interval (grey shaded area) for $b\mapsto \widehat\EE_b(V,T) $ with $T=20000$ under the condition $\widehat\EE_b(V,T)<\infty$. Also included are $b \mapsto \EE_b(V)$ (black dots), and the theoretical finite maximum  $(T-1)/2$ (horizontal red dashed line).}
    \label{fig:bias_stat}
\end{figure}
We then consider $\widehat\EE_b(V,T)$ for $T=20000$, $b\in[2.5,4.5]$, the resulting mean and upper/lower $10\%$ quantiles are plotted in \Cref{fig:bias_stat}. Note that by \eqref{eq:statmean} we have that $\widehat\EE_b(V,T)<\infty$ implies $\widehat\EE_b(V,T)\leq (T-1)/2$, which provides an upper bound on the finite values the estimator can take. As a specific consequence, if $\EE_b(V)> (T-1)/2$ then $\widehat\EE_b(V,T)<\EE_b(V,T)$ and there must necessarily exist a negative bias. Similarly, in \Cref{fig:bias_stat}, we see that the conditional mean of $\widehat\EE_b(V,T)$ deviates from $\EE_b(V)$ only when $\EE_b(V)$ approaches $(T-1)/2$.

\subsection{Non-Stationary Case}

In the case where $V$ is not stationary we again consider the AR1 process defined in \eqref{eq:AR1_V}, but scale it by a periodic function. Specifically 
\begin{align*}
    V_t&=Z_t S_t,\\
    Z_t&=\sqrt{1-\rho^2}\sum_{s=-\infty}^t \rho^{t-s}X_s,\\
    S_t&=\frac{1}{2}{\left(1+\sin\left(\frac{t-\lfloor t/1000\rfloor }{1000} \,\pi \right)^2\right)}
\end{align*}
for all $t\in\Z$. Here $\rho=70\%$ and $\{X_t\}_{t\in\Z}$ is an independent sequence of standard normal random variables. The function $S$ starts at $1$ for $t=0$, and then drops down to $0.5$ for $t=500$, finally returning to $1$ for $t=1000$. We also have $S$ periodic with period $1000$, thereby imitating a seasonal scale for $V$.

We may then apply employ \Cref{prop:nonstatestimator} in order to estimate $\EE_b$ by
$$\EE_b(V)=\E\left[\tau_{(b,\infty)}(V,0)\right]=\E\left[\tau_{\alpha}(Z,0)\right],$$
where $\alpha_t=[b/S_t,\infty]$.

\begin{figure}[h]
    \centering
    \includegraphics[width=0.75\textwidth]{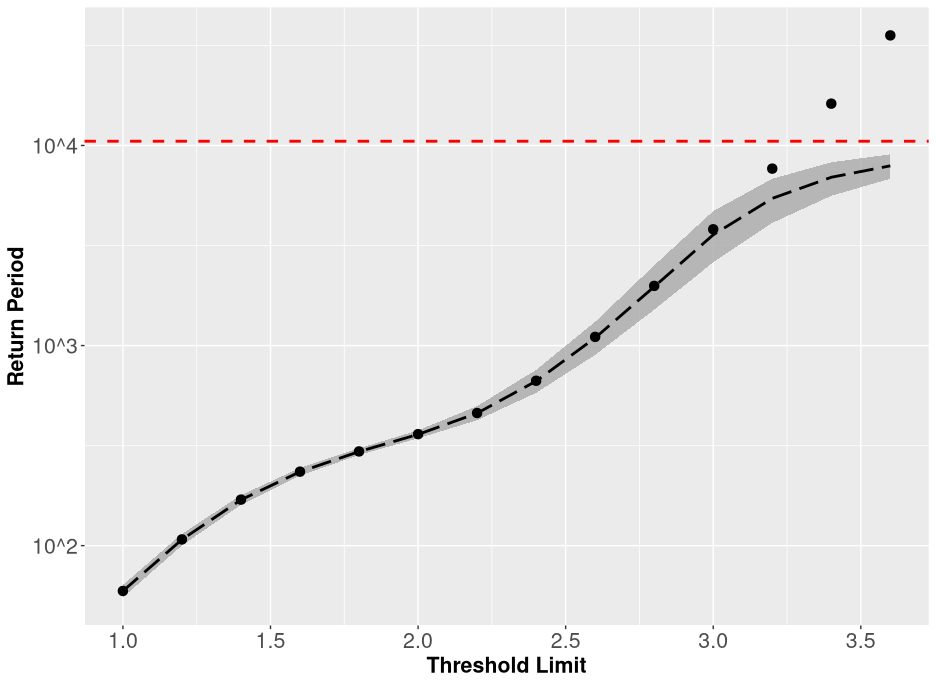}
    \caption{Plot of conditional mean (black dashed line) and $80\%$ prediction interval (grey shaded area) for $b\mapsto \widehat\EE_b(V,T) $ with $T=20000$ under the condition $\widehat\EE_b(V,T)<\infty$. Also included are $b \mapsto \EE_b(V)$ (black dots), and the theoretical finite maximum  $(T+Y-2)/2$ (horizontal red dashed line).}
    \label{fig:bias_nonstat}
\end{figure}

We then consider $\widehat\EE_b(V,T)$ for $T=20000$, $b\in[1,3.5]$, the resulting conditional mean and conditional upper/lower $10\%$ quantiles are plotted in \Cref{fig:bias_nonstat}. We have also included the line corresponding to $(T+Y-2)/2$, where $Y=1000$ is the period of $S$, representing the length of a \textit{year}. In this case it can be computed that, conditional on $\widehat\EE_b(V,T) < \infty$, we have $\widehat\EE_b(V,T) \leq (T+Y-2)/2 $ for any $b$.

In examining \Cref{fig:bias_nonstat}, we see many similarities with \Cref{fig:bias_stat}. Specifically the estimator is unbiased until $\EE_b$ approaches the maximum limit of $(T+Y-2)/2$. 

It is common to ignore the effects of seasonality, assuming that the seasonal effects will in some sense average out over time. To examine the effect of the inclusion of seasonality in the estimation of $\EE_b(V)$ we additionally compute our estimator under the assumption that $V$ is stationary by applying \eqref{eq:statmean}. 

For any $b\in\R$, we recall the definition  $ \Theta\triangleq\{\theta\in\N: \theta<T,\, V_\theta > b\} $. We denote by $ \theta_i $, $ i\leq|\Theta| $ the $ i $'th element of $ \Theta $, and define $ \theta_{|\Theta|+1}=T+\theta_1 $. We lastly introduce the shorthand $\Delta\theta_i=\theta_{i+1}-\theta_i$. With this we denote
\begin{equation}\label{eq:statest_in_nonstat}
    \widetilde\EE_b(V,T) =  \sum_{i=1}^{|\Theta|}
            \frac{\Delta\theta_i(\Delta\theta_i-1)}{2 T},
\end{equation}
unless $\Theta=\emptyset$ in which case we set $\widetilde\EE_b(V,T) = \infty $.

\begin{figure}[h]
    \centering
    \includegraphics[width=0.75\textwidth]{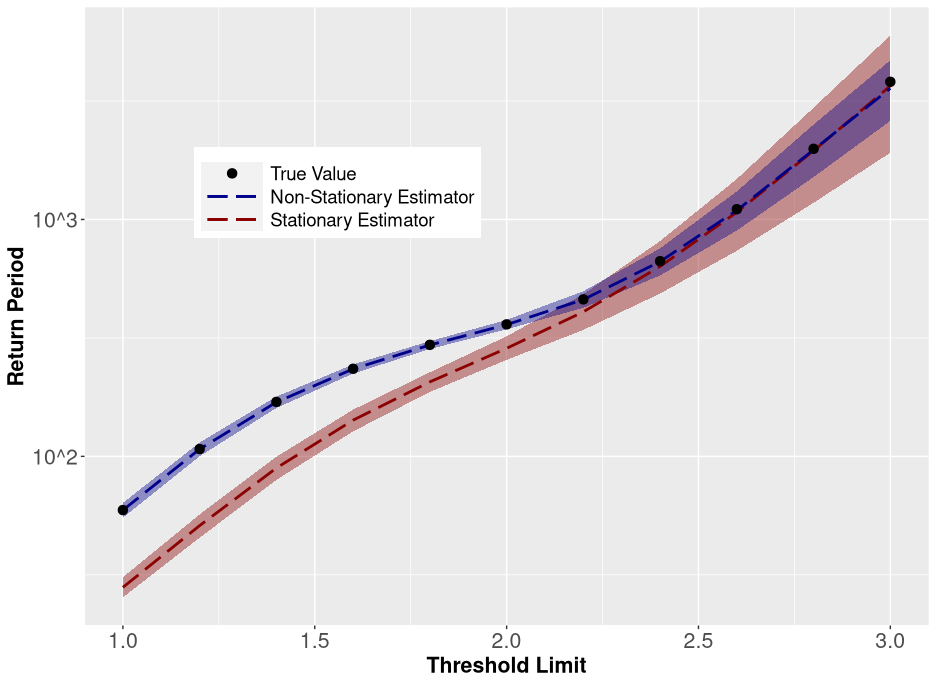}
    \caption{Plot of conditional mean and $80\%$ prediction interval for $b\mapsto \widehat\EE_b(V,T) $ and $b\mapsto \widetilde\EE_b(V,T) $  for $T=20000$ (all conditional on having finite values),  along with $b \mapsto \EE_b(V)$.} 
    \label{fig:bias_stat_nonstat_compare}
\end{figure}

Exceedances are less likely to occur for $t$ close to $T\Z$, as these times correspond to the lowest values of $S$.  Therefore, there is a lower chance of exceedance for low values of $t$, leading the non-stationary estimator to usually return higher values than the stationary for low values of $b$, as apparent in \Cref{fig:bias_stat_nonstat_compare}. As we see, this bias disappears as we approach the point where $\EE(V)=Y$. Generalising this effect we can easily argue that a stationary model for $V$ is sufficient as long as we consider only larger return periods. This will also have the added benefit of avoiding potential issues with overfitting.

Despite this, we may also note from \Cref{fig:bias_stat_nonstat_compare} that the prediction interval for the non-stationary estimator is a lot tighter than the non-stationary one, implying more precise estimates. This can be explained intuitively, by the fact that using the additional knowledge of the precise seasonal behaviour of $V$ should yield better estimates, even for larger return periods. As such, this encourages the inclusion of seasonal behaviour whenever it can be accurately estimated.

As a final note, it is worth mentioning the case where the non-stationarity of $V$ is not purely seasonal, but also includes a long-term trend. For example, in the case of wave heights, there has been a significant increase over the years, see e.g.\ \cite{Kushnir,vanemBayes}. In cases with long-term trends, the non-stationarity of $V$ might have a big impact on the tail behaviour of exceedance times. As such we would have no guarantee that a stationary estimate approaches the true return period for high thresholds as we see in the purely seasonal case. The inclusion of non-stationarity in the modelling of $V$ is therefore particularly important in such cases.

\section{Comparison of Methods and Models for Empirical Data}\label{sec:empassessment}
After seeing how the estimator performs on synthetic data we can move on to empirical analyses. The goal will be to compare the average exceedance times of four different models with empirical estimates. The data considered for this example will be ERA5 reanalysis data \cite{ERA5}. Specifically, we will use hourly wind speeds in meters per second (m/s) at 100 meters elevation in an area located at 49\degree N, -8\degree W. The data is for every hour over the 41-year period $ 1980$-$2020 $, yielding a total of 359\,424 individual data points.

We again consider the quantities
\begin{equation*}
    \EE_b(V) \triangleq \E[\tau_{(b,\infty)}(V,0)], \quad 
    \widehat\EE_b(V,T) \triangleq \frac{1}{T}\sum_{t=0}^{T-1}\widehat\tau_{(b,\infty)}(V,t).
\end{equation*}
The idea is to use $\widehat\EE_b(V,T)$ and its corresponding confidence interval as an estimate for $\EE_b(V)$. Using this estimate as an approximate ground truth, we can evaluate and compare the performance of different models of $V$. This is done by computing $\EE_b(V^*)$, for some $V^*$ taken from a selection of models for $V$.

\subsection{Stationary Models}

For our initial analysis we will ignore the seasonal trends present in the data and treat $V$ as an ergodic process {satisfying the conditions of our central limit theorem (\Cref{thm:feasibleCLT}). This allows us to compute $b\mapsto\widehat\EE_b(V,T)$, $T=359424$, along with an approximate $95\%$ confidence interval for $\EE_b(V)$.}



In providing models for $V$ we will here focus on two aspects. Firstly, we require the marginal distributions of $V$ in the form of the cumulative distribution function $F_V(v)=\PP(V_t\leq v)$. Note that as a consequence of ergodicity, this function is independent of $t$. Secondly, we will consider the autocorrelation structure $\rho(V,t)=\text{cor}(V_0,V_t)$. To compute $F_V$ we assume that
$$F_V(v)=1-\text{exp}\left(-\left(\frac{v}{\lambda}\right)^\alpha\right),$$
i.e.\ that the marginals of $V$ follow a Weibull distribution with scale $\lambda=11.05$ and shape $\alpha=2.19$. Similarly, we estimate the autocorrelation structure of $V$ as
\begin{equation}\label{eq:rho_stat}
\rho(V,t)=\left(1+\kappa\left(\frac{t}{\zeta}\right)^\eta\right)^{-\frac{1}{\eta\kappa}},
\end{equation}
with $\zeta=10.23, \, \eta=1.63, \, \kappa=1.38$. Note that several parametric models for $\rho(V,\cdot)$ were considered, see e.g.\ \cite{papalexiou2018unified} for a list of suggested models. However, \eqref{eq:rho_stat} was found to provide the best fit for our dataset. We also remark that the same parametric model was used for the autocorrelation structure of significant wave heights in \cite{vanem2023analysing}.

{We will in particular consider a class of models known as Gaussian copula processes. Let $\Phi$ be the cumulative distribution function of a standard normal Gaussian. We then have that $Z_t =\Phi^{-1}F_V(V_t)$ has a standard normal distribution for any $t$. It is important to note that this does not in general make $Z$ a Gaussian process. However, modelling it as such allows for simple models that still capture the correct marginal distributions and autocorrelation structure of $V$. This approach was suggested for the modelling of environmental factors in \cite{papalexiou2018unified}, and has been applied to compute return periods in the context of environmental contours \cite{vanem2023analysing,vanem2023analyzing_nonstat}.}

Following \cite{papalexiou2018unified,vanem2023analysing} we define $\TT$ as the transformation $\TT(\rho(Z,t))=\rho(V,t)$, which can be computed as
$$\TT(\rho)=\frac{\int_\R\int_\R F_V^{-1}\Phi(x)F_V^{-1}\Phi(y)\phi(\rho;x,y)  dxdy -\mu_V^2  }{\sigma_V^2}, $$
where $\phi(\rho;\cdot,\cdot)$ is the density function of a bivariate normal distribution with correlation $\rho$. Additionally, $\mu_V=9.80$ and $\sigma_V=4.71$ are the mean and standard deviation, respectively, of $V$. We then aim to use the relation $\TT(\rho(Z,t))=\rho(V,t)$, which is done in e.g.\ \cite{vanem2023analysing,papalexiou2018unified}, by using an accurate parametric approximation of $\TT$. In our case it was found that 
\begin{equation}\label{eq:cor_transform}
\TT(\rho)\approx \frac{(1+\xi\rho)^\upsilon-1}{(1+\xi)^\upsilon-1}, 
\quad
\TT^{-1}(\rho) \approx \frac{(1+\rho(1+\xi)^\upsilon-\rho)^{1/\upsilon}-1}{\xi},
\end{equation}
for $\xi=0.065, \, \upsilon=0.373$. Note that under the assumption that $Z$ is a strictly stationary standardised Gaussian process, it will be uniquely determined by $\rho(Z,\cdot)$.

The first model we consider, labeled $V^{(\text{IID})}$, will be the simplest one. We define $\{V^{(\text{IID})}_t\}_{t\in\Z}$ to be an i.i.d.\ sequence of random variables such that $V^{(\text{IID})}_t=V_t$ in law for all $t\in\Z$. Consequently we can compute the average exceedance times of $V^{(\text{IID})}$ by
$$\EE_b\left(V^{(\text{IID})}\right)=\frac{F_V(b)}{1-F_V(b)}=\text{exp}\left(\left(\frac{b}{\lambda}\right)^\alpha\right)-1.$$
Note that we could equivalently consider $Z^{(\text{IID})}=\Phi^{-1}F_V\left(V^{(\text{IID})}\right)$, where $\{Z^{(\text{IID})}_t\}_{t\in\Z}$ is a sequence of i.i.d.\ standard normal variables.

For $V^{(\GG)}$ we consider a general standard Gaussian process $Z^{(\GG)}$ such that $V^{(\GG)}=F_V^{-1}\Phi(Z^{(\GG)})$ has the same autocorrelation structure as $V$, i.e.\ $\rho(V^{(\GG)},t)=\rho(V,t)$.  This allows us to compute $\rho\left(Z^{(\GG)},t\right)=\TT^{-1}\rho(V,t)$ by combining \eqref{eq:rho_stat} and \eqref{eq:cor_transform}. Once $\rho\left(Z^{(\GG)},t\right)$ is computed we can compute $\EE_b\left(V^{(\GG)}\right)$ by simulation of $Z^{(\GG)}$.

{We further define $V^{(\text{AR1})}=F_V^{-1}\Phi Z^{(\text{AR1})}$. Here, $Z^{(\text{AR1})}$ is the unique standardised Gaussian process with $\rho\left(Z^{(\text{AR1})},t\right)=\exp(-\theta t)$, $\theta=2.7\%$ for any $t\in\N$. This makes $Z^{(\text{AR1})}$ a Gaussian AR1 process. Note that $\theta=2.7\%$ was chosen to minimise $\sum_{n=1}^{100} (\rho(V,n)-\rho(V^{(\text{AR1})},n))^2$. This simpler process can have several advantages over $Z^{(\GG)}$, such as being Markovian, which makes for more efficient simulation. Additionally, in \cite{huseby2023AR1}, this model was leveraged to enable an importance sampling technique for the construction of environmental contours.}

\begin{rem}
    It is common to calibrate AR1 models by considering only $\exp(-\theta)=\rho\left(Z^{(\text{AR1})},1\right)=\TT^{-1}\rho\left(V,1\right)$ which only requires $\rho\left(V,1\right)$. However, we can see from \eqref{eq:rho_stat} and \eqref{eq:cor_transform} that $\TT^{-1}\rho\left(V,t\right) < (\TT^{-1}\rho\left(V,1\right))^t$ for all $1\leq t \leq 100$, which means that calibration based only on $\rho\left(V,1\right)$ could significantly overestimate the autocorrelation. For reference, $\TT^{-1}\rho\left(V,1\right)=\exp(-1.4\%)$ which would nearly halve the value of $\theta$ compared to the chosen method.
\end{rem}

We can also consider a model in continuous time. Specifically, we define $V^{(\text{OU})}(t)=F_V^{-1}\Phi(Z^{(\text{OU})}_t)$ with
$$Z^{(\text{OU})}(t)=\sqrt{2\theta}\int_{-\infty}^t\exp\left(\theta(t-s)\right)dW_s,$$
for any $t\in\R$ where $W$ is a standard Wiener process. This process is known as a standardised Ornstein-Uhlenbeck process and serves as a continuous interpolation of an AR1 process. In particular, the discrete process $t\mapsto Z^{(\text{OU})}(t), t\in\Z$ equals $Z^{(\text{AR1})}$ in distribution. We may also extend $\EE_b(\cdot)$ to our continuous-time process by $\EE_b(V^{(\text{OU})})=\E\left[\inf\{t\in\R,\,t\geq0:V^{(\text{OU})}(t)>b\}\right]$. An explicit expression for this value (see e.g.\ \cite{OUmom}) is given by
\begin{align*}
    \EE_b\left(V^{(\text{OU})}\right)&= \E\left[\left( \phi\left(\Phi^{-1}F_V(b)\right)-\phi\left(Z^{(\text{OU})}_0\right) \right)\1(V^{(\text{OU})}_0<b)\right],\\
    \phi(x) &= \frac{1}{2\theta}\sum_{i=1}^{\infty} \frac{ (\sqrt{2}x)^{i}}{i!}\Gamma\left(\frac{i}{2}\right). 
\end{align*}
We also have a more convenient formulation in
\begin{equation*}
    \EE_b\left(V^{(\text{OU})}\right)
    =\E\left[{\frac{\sqrt{2\pi}}{\theta}}\int_{Z^{(\text{OU})}_0 \land \Phi^{-1}F_V(b)}^{\Phi^{-1}F_V(b)}\Phi(t)e^{t^2/2}dt\right],
\end{equation*}
which can be readily computed by Monte-Carlo simulation of $Z^{(\text{OU})}_0$, which has a standard normal distribution. Alternatively, we have the deterministic asymptotic approximation $\EE_b\left(V^{(\text{OU})}\right) \approx \sqrt{2\pi}\exp((\Phi^{-1}F_V(b))^2/2)/(\theta\Phi^{-1}F_V(b)) $, for high values of $b$.

\begin{figure}[h]
    \centering
    \includegraphics[width=0.95\textwidth]{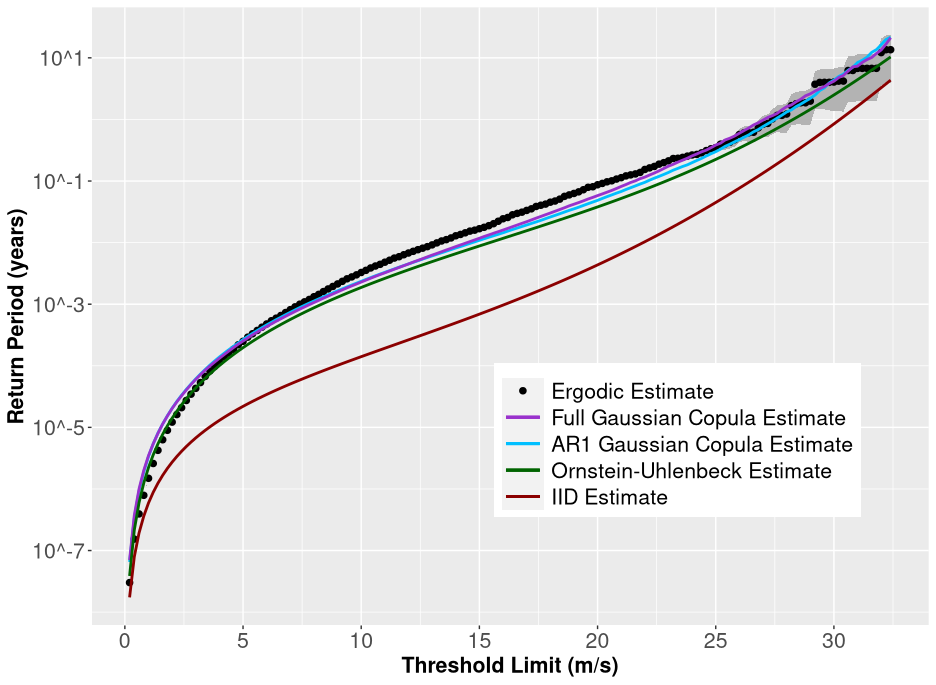}
    \caption{Comparison of $b \mapsto \widehat\EE_b(V,T)$ and $b \mapsto \EE_b(V^{*})$ for $V^*=V^{(\GG)}$, $V^{(\text{AR1})}$, $V^{(\text{OU})}$, and $V^{(\text{IID})}$}
    \label{fig:compare_stationary}
\end{figure}

The return periods of the different models, $b\mapsto \E_b(V^{(*)})$, $V^{(*)}=V^{(\GG)}$, $V^{(\text{AR1})}$, $V^{(\text{OU})}$, and $V^{(\text{IID})}$, $b\in [0.4,32.4]$ (m/s), are all plotted in \Cref{fig:compare_stationary}, along with the estimated return periods $b\mapsto\widehat\EE_b(V,T)$. What we see is that most of the models perform quite well in estimating the true return period $\EE_b(V)$, at least in regions where we can accurately estimate it by $\widehat\EE_b(V,T)$.

It appears that the best fit is given by $V^{(\GG)}$, which captures the full autocorrelation structure of $V$. Despite this we may note that $V^{(\text{AR1})}$, which uses a simpler and more approximate autocorrelation, gives almost identical results. This may indicate that the behaviour of $\tau_\beta(V,\cdot)$ can be accurately captured by Markovian processes, which opens up for more efficient methods that take advantage of this Markovianity.

Our other Markovian process, $V^{(\text{OU})}$, performs slightly worse. Since $V^{(\text{OU})}$ is a continuous-time interpolation of $V^{(\text{AR1})}$ we must have $\EE_b(V^{(\text{OU})}) < \EE_b(V^{(\text{AR1})})$, which inherently makes $V^{(\text{OU})}$ a more conservative model. {It is also possible to consider $\EE_b(V^{(\text{OU})})$ as an approximation of $\EE_b(V^{(\text{AR1})})$ in the sense that their ratio approaches 1 for $\theta \to 0$. Consequently, we can use the explicit formula for $\EE_b(V^{(\text{OU})})$ to efficiently compute a conservative estimate of $\EE_b(V^{(\text{AR1})})$, and hence $\EE_b(V)$.}

We have also included $V^{(\text{IID})}$ as it is a commonly used simple and conservative model for the purposes of computing return periods. We here see that this model significantly underestimates the return period. In fact, one could consider this phenomenon theoretically by applying Slepian's lemma. This result states that if $\rho(X,\cdot)\leq \rho(Y,\cdot)$ for any two stationary standard normal Gaussian processes $X,Y$ then
$$    \PP\left(\max_{0\leq t \leq t^*}X_t \leq x \right) 
\leq \PP\left(\max_{0\leq t \leq t^*}Y_t \leq x\right),$$
for all $t^*\in\N$. Since $\rho(Z^{(\text{IID})},t)=0$ for any $t>0$, and $0<\rho(Z^{(\GG)},\cdot),\rho(Z^{(\text{AR1})},\cdot)$, we get $\EE_b(V^{(\text{IID})})\leq \EE_b(V^{(\GG)}),\EE_b(Z^{(\text{AR1})})$ for all $b$. Comparison is more difficult for $Z^{(\text{IID})}$ and $Z^{(\text{OU})}$, as $Z^{(\text{OU})}$ is defined in continuous time. However, in \cite{sande2024convex}, it was shown that there exists some threshold $b^*$ such that $\EE_b(V^{(\text{IID})}) > \EE_b(V^{(\text{OU})})$ for all $b>b^*$, while we generally see $\EE_b(V^{(\text{IID})}) < \EE_b(V^{(\text{OU})})$ for all $b<b^*$. For the calibrated value of $\theta$ we have $b^*\approx $ 44.6 m/s, where we have $\EE_{b^*}(V^{(\text{IID})}) \approx \EE_{b^*}(V^{(\text{OU})})\approx $ 1.7e7 years.

\subsection{Non-Stationary Models}
As mentioned, there is a significant seasonal effect present in our dataset. In order to properly include this seasonality in our models we consider the following shape of the marginal distribution of $V$.
$$V_t\sim W(\lambda_t,k_t),$$
i.e.\ $V_t$ is Weibull distributed with scale $\lambda_t$ and shape $k_t$. We further assume that $k$ and $\lambda$ are periodic functions with a period of one year. We may then model the cumulative distribution function of $V_t$ by 
$$F_t(x)=1-\exp\left(-\left(\frac{x}{\lambda_t}\right)^{k_t}\right). $$

This implies that $ U_t\triangleq F_t(V_t)$ is constant in distribution. If we further assume that $U$ is ergodic, we may apply \Cref{prop:nonstatestimator} to compute $\E[\tau_\beta(V,0)]$ by $\E[\tau_{\alpha}(U,0)]$, with $\alpha_s=F_{s}(\beta_s)$.

In calibrating $\lambda$ and $k$ we employ a similar method to what was used in \cite{sande2024convex}. To calibrate $k$ we invert the equality
\begin{equation*}
    \frac{\Gamma(1+1/k(t))^2}{\Gamma(1+2/k(t))}= \frac{\E[V_t]^2}{\E[V_t^{2}]},
\end{equation*}
where $\Gamma$ is the gamma function, and the expectations are computed using spline regression. We can then estimate $\lambda$ by
\begin{equation*}
    \lambda_t = \E\left[\frac{V_t}{\Gamma(1/k(t))}\right]
\end{equation*}

\begin{figure}[h]
	\centering
	\begin{subfigure}{.5\textwidth}
		\centering
		\includegraphics[width=.95\linewidth]{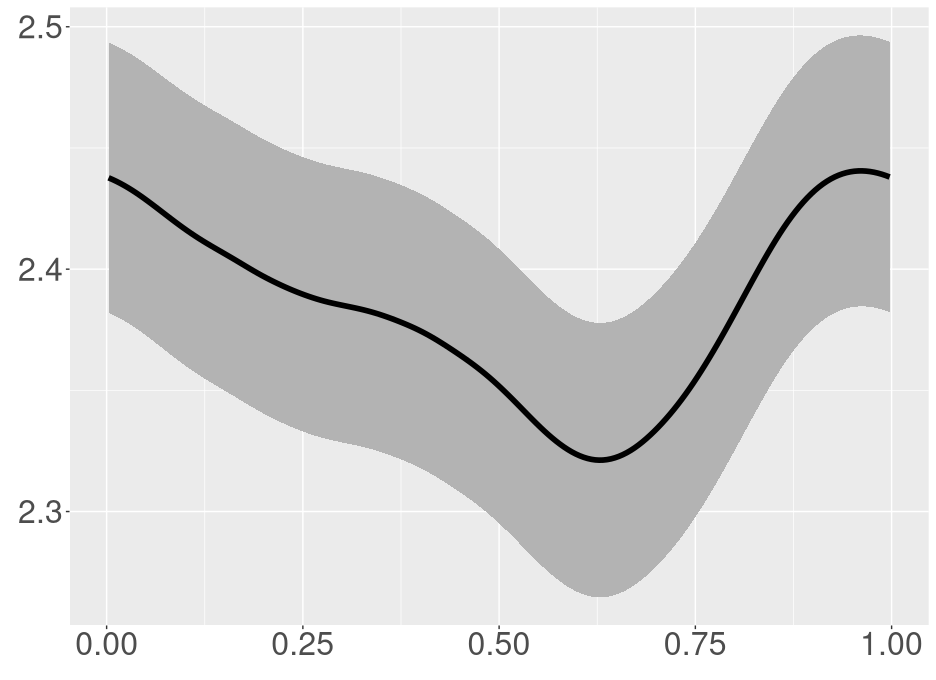}
		\caption{Values of $ k $}
		\label{fig:ks}
	\end{subfigure}%
	\begin{subfigure}{.5\textwidth}
		\centering
		\includegraphics[width=.95\linewidth]{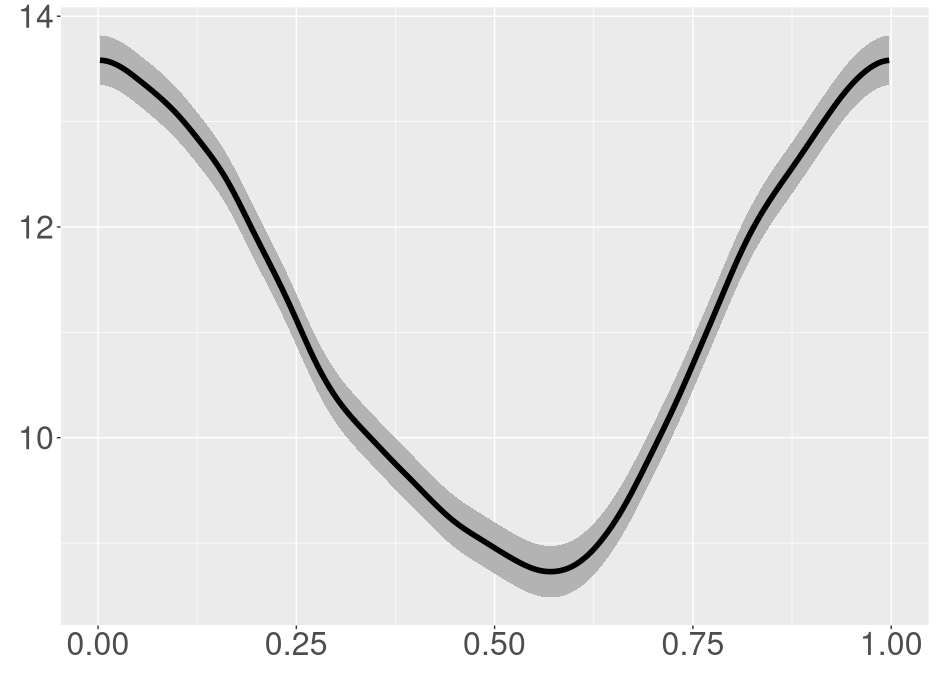}
		\caption{Values of $ \lambda $}
		\label{fig:lams}
	\end{subfigure}
	\caption{Non-parametric functions for $ F_t $, with 95\% confidence intervals.}
	\label{fig:Hfuncs}
\end{figure}

We again consider the non-stationary analogue of the i.i.d.\ model where $\{V^{\text{IS}}_t\}_{t\in\Z}$ is an independent sequence of random variables such that $V^{\text{IS}}_t=V_t$ in distribution.

For the second non-stationary model we consider $t\mapsto F_t^{-1}\Phi(Z^{(\GG 2)}_t)$ where $Z^{(\GG 2)}$ is a stationary standardised Gaussian process. We parametrise the autocorrelation of $Z^{(\GG 2)}$ by
$$\rho(Z^{(\GG 2)},t)=\left(1+\kappa_2\left(\frac{t}{\zeta_2}\right)^{\eta_2}\right)^{-\frac{1}{\eta_2\kappa_2}},$$
where $\zeta_2=10.65, \, \eta_2=1.56, \, \kappa_2=0.83$. These parameters are chosen to minimise
\begin{equation*}
    \sum_{n=1}^{100} \left(\rho(Z,n)-\rho\left(Z^{(\GG 2)},n\right)\right)^2,
\end{equation*}
where $Z_t=\Phi^{-1}F_t(V_t)$.

Lastly, we consider a non-Gaussian copula process $V^{(\mathbf{t})}_t= F_t^{-1}\mathbf{t}_{\nu}(Y^{(\mathbf{t})}_t)$, where $\mathbf{t}_\nu$ is the student's t cumulative distribution function with degrees of freedom $\nu$. Here, $Y^{(\mathbf t)}$ is the unique stationary Markovian process such that the joint distribution of $(Y^{(\mathbf{t})}_0,Y^{(\mathbf{t})}_1)$ is a bivariate t-distribution with degrees of freedom $\nu=13.4$ and correlation coefficient $\rho_\mathbf{t}=96.4\%$. 

{In order to compute the specified values for $\nu$ and $\rho_\mathbf{t}$ we consider the following.}
The aim is to make $U^{(\mathbf{t})} = \mathbf{t}_{\nu}(Y^{(\mathbf{t})})$ as similar to $U$ as possible. We first define the tail dependence coefficient for a stationary process $X$ and threshold $b$ as
\begin{equation}\label{eq:taildepdef}
    \lambda(X,b) \triangleq \PP\left(X_0>b \big| X_{-1} > b \right).
\end{equation}
As seen in e.g.\ \cite{demarta2005t}, we have 
\begin{equation}
    \lim_{b\to 1}\lambda(U^{(\mathbf{t})},b)= 2\mathbf{t}_{\nu+1}\left(-\sqrt{\frac{(\nu+1)(1-\rho_\mathbf{t})}{1+\rho_\mathbf{t}}}\right)
\end{equation}
The parameters $(\nu,\rho_\mathbf{t})$ were then chosen to minimise 
$$\sum_{n=1}^{100} (\rho(U,n)-\rho(U^{(\mathbf{t})},n))^2$$
under the constraint of 
$\lim_{b\to 1}\lambda(U^{(\mathbf{t})},b)=\lambda(U,\text{1-1e-4}) = 61.5\%$.

\begin{figure}[h]
    \centering
    \includegraphics[width=0.95\textwidth]{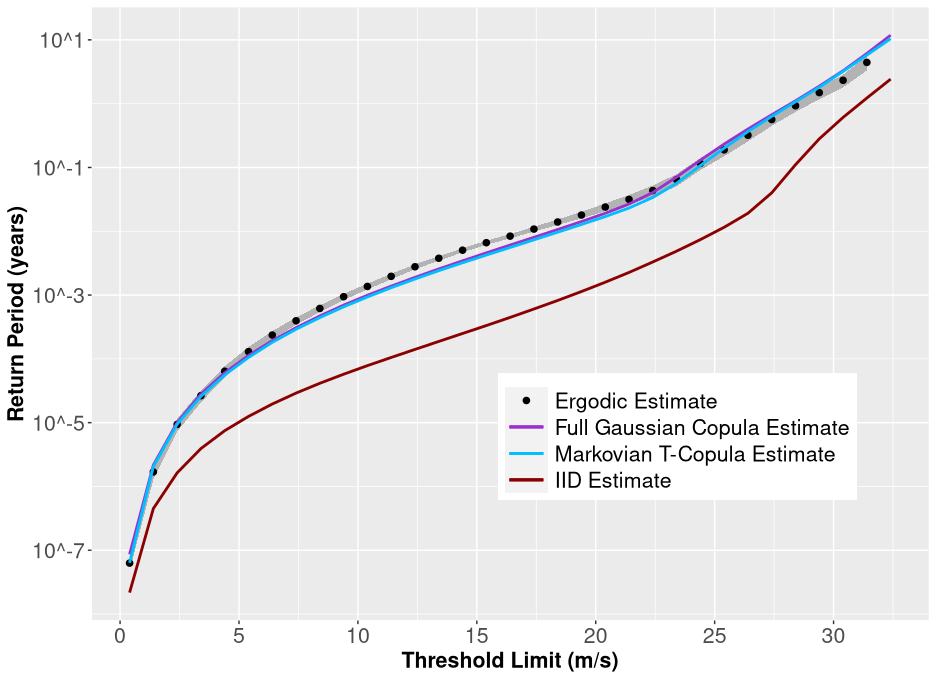}
    \caption{Comparison of $b \mapsto \widehat\EE_b(V,T)$ and $b \mapsto \EE_b(V^{*})$ for $V^*=V^{(\GG 2)}$, $V^{(\text{T})}$, and $V^{(\text{IS})}$}
    \label{fig:compare_nonstat}
\end{figure}

In \Cref{fig:compare_nonstat}, we see very similar results as in \Cref{fig:compare_stationary}. The model based on an independent sequence, $\EE_b(V^{(\text{IS})})$, greatly underestimates $\EE_b(V)$, in contrast to $\EE_b(V^{({\GG 2})})$ and $\EE_b(V^{({\mathbf{t}})})$, which more accurately captures the correlation structure of $V$. This indicates that the extreme behaviour of $V$ can be accurately captured by our modelling copula processes even in our non-stationary setting.

It is worth mentioning that while $V^{({\GG 2})}$ and $V^{(\mathbf{t})}$ provide nearly identical return periods, they have quite different asymptotic properties. To examine these asymptotics we first denote $U^{(\mathbf{t})} = \mathbf{t}_{\nu} (Y^{(\mathbf{t})})$ and $U^{(\GG 2)} = \Phi (Z^{(\GG 2)})$ in order to reconcile their marginal distributions. We previously mentioned that $\lim_{b\to 1}\lambda(U^{(\mathbf{t})},b)=61.5\%$, and since $Z^{(\GG 2)}$ is a Gaussian process we further have that $\lim_{b\to 1}\lambda(U^{(\GG 2)},b)=0$. This asymptotic tail dependence is directly related to another important concept, namely the \textit{up-crossing rate}, defined for any stationary process, $X$, as
\begin{equation}\label{eq:upcrossdef}
    \mu(X,b) \triangleq \PP\left(X_0>b ,\, X_{-1} \leq b\right).
\end{equation}
To see how \eqref{eq:taildepdef} and \eqref{eq:upcrossdef} are related we note the following. For any stationary process $X$ we have
\begin{align*}
    \PP(X_0 > b))
    =& \PP\left(X_0>b ,\, X_{-1} \leq b\right)\\
    &+\PP\left(X_0>b \big| X_{-1} \leq b \right)\PP\left(X_{-1}>b\right)\\
    =&\mu(X,b) + \lambda(X,b)\PP(X_0 > b)),
\end{align*}
which implies 
\begin{equation}\label{eq:tail_to_upcross}
    1-\lambda(X,b)=\frac{\mu(X,b)}{\PP(X_0 > b))}.
\end{equation}
A particular consequence of this is
\begin{equation}\label{eq:ratio_exceedrate}
    \lim_{b\to1}\frac{\mu(U^{(\mathbf{t})},b)}{\mu(U^{(\GG 2)},b)}
    =\lim_{b\to1}\frac{\mu(U^{(\mathbf{t})},b)}{\mu(U^{(\text{IS})},b)}
    = 1-\lim_{b\to1}\lambda(U^{(\mathbf{t})},b)
    = 38.5\%
\end{equation}
A higher upcrossing rate usually implies a lower average time until next exceedance. As such we can expect that $ \EE_b(V^{({\GG 2})}) / \EE_b(V^{({\mathbf{t}})})$ will behave similar to \eqref{eq:ratio_exceedrate} for sufficiently high values of $b$. 

As such, we have two models with significantly different asymptotic properties, that mostly agree in our chosen regime of moderate values of $b$. This serves as a reminder that model performance in our chosen regime is insufficient to guarantee that asymptotic properties of the true model are accurately captured. 

We may finally note that since $\lambda(X,b)\to 0$ for any ergodic Gaussian process $X$, we have $\lambda(X,b) \leq \lambda(Z,b)$ for sufficiently large $b$, implying $\mu(X,b) \geq \mu(Z,b)$. Since $X$ will exhibit more frequent upcrossings than $Z$ we can usually expect that $\EE_b(X)<\EE_b(Z)$. Consequently, we can consider the use of Gaussian copula models as an asymptotically conservative estimate of $\EE_b(V)$. 

In summary, Gaussian copula models can accurately capture $\EE_b(V)$ in our chosen regime, while providing a conservative estimate in the asymptotic regime, where data is sparse to non-existent.

\begin{rem}
    In order to expand on the implications of the difference in asymptotic properties for $U^{(\mathbf{t})}$ and $U^{(\GG 2)}$ we consider two concepts closely related to tail dependence, the extremal index and cluster sizes.

    To define the extremal index of a stationary process, $\theta(X)$, we assume the following. If, for any $v\in \R^+$ and sequence $\{u_n^{(v)}\}_{n\in\N}\subset\R$ such that 
    $$\lim_{n\to\infty}n\PP(X_0 > u_n^{(v)}) = v$$
    then 
    $$\lim_{n\to\infty}\PP(\max_{t=1,\cdots,n}X_t \leq u_n^{(v)}) = e^{-\theta(X) v},$$
    if the extremal index exists. For a more comprehensive definition and discussion of the extremal index, including conditions for existence, we refer to \cite{beirlant2006statistics}.

    By cluster size, we refer to the average time spent above a threshold $b$ after an upcrossing. Specifically, for any stationary process $X$, we define 
    $$C(X,b)=\E[\min \{i: X_i\leq b \} | X_0>b,X_{-1}\leq b].$$

    In \cite{MartaFerreira2015} they discuss estimating the extremal index by means of the tail dependence. For example, if $V$ satisfies a specific technical condition ($D^2(u_n^{(v)})$, as introduced in \cite{chernick1991calculating}), it was proved in \cite{ferreira2012extremal} that $\theta(V)=1-\lim_{b\to\infty}\lambda(V,b)$.

    Similarly, in \cite{meiss1997average}, the relation $C(X,b)^{-1}=\mu(X,b)/\PP(X_0>b)$ was shown for a general ergodic process $X$. By \eqref{eq:tail_to_upcross} we therefore have $C(X,b)^{-1}=1-\lambda(X,b)$. Under $D^2(u_n^{(v)})$ we would lastly get the relation $\theta(X)^{-1}=\lim_{b\to\infty}C(X,b)$.
\end{rem}


\appendix
\section{Technical Results and Proofs}\label{sec:resultsandproofs}


	
	%
	%
		%


In this section we aim to provide a rigorous proof of convergence for our estimator
$$\E[f(\tau_\beta(V,0))]\approx \sum_{t=0}^{T-1}f(\widehat\tau_\beta(V,t))$$
where $f$ is any function such that $\E[|f(\tau_\beta(V,0))|]<\infty$ and
\begin{equation*}
\widehat\tau_\beta(T;V,t) =
\begin{cases}%
  \tau_\beta(V,t) & \tau_\beta(V,t) < T-t\\
    T-t+\Lambda(T;V,t) & \tau_\beta(V,t) \geq T-t.
\end{cases}
\end{equation*}
%

Note that we, for the sake of generality, no longer assume the specific shape of $\Lambda$ given in \eqref{eq:Lamdadef}. 

To guarantee convergence of our estimator we first introduce the concept of \textit{transit times}.
\begin{defi}
    Assume $\widetilde\beta\subset \R^N$. For any process $V$ and time $T$ we may define 
    $$ \theta_{\widetilde\beta}^-(V,T) \triangleq \max(t \leq T:V_t \in \widetilde{\beta}), \quad \theta_{\widetilde\beta}^+(V,T) \triangleq \min(t \geq T:V_t \in \widetilde{\beta}).$$ 
    Note the ordering $\theta_{\widetilde\beta}^-(V,T)\leq T \leq \theta_{\widetilde\beta}^+(V,T)$.
    We can then further define the \textit{transit time} at time $T\in\Z$ by, 
    $$ \tau^\pm_{\widetilde\beta}(V,T) \triangleq \max(\theta_{\widetilde\beta}^+(V,T)-\theta_{\widetilde\beta}^-(V,T) - 1,0).$$
\end{defi}


Using this concept we can bound the difference between sums of $f(\tau_\beta(V,t))$ and sums of $f(\widehat\tau_\beta(V,t))$

\begin{lemm}\label{lem:bound_by_P}
	Assume the threshold function, $ \beta $, satisfies $ \beta_s\supseteq\widetilde{\beta} $ for some measurable $\widetilde\beta\subseteq\R^N$ and all $ s\geq 0 $. Further fix a function $ f:\N\mapsto\R $ and assume that
	$$ |f(t)| \leq C_1 + C_2t^{p}, $$
	for some $ p,\, C_1,\,C_2\in[0,\infty) $. 
 
	Under these assumptions we have, for any $q\in[0,\infty)$ that
	$$ \sum_{t=0}^{T-1}\big|f(\tau_\beta(V,t)) - f(\widehat{\tau}_\beta(T;V,t))\big|^q
    \leq \Pp_{p,q}(T,\overline \Lambda _T).$$
    Here
        \begin{equation*}
            \Pp_{p,q}(T,C)=
            { \tau^\pm_{\widetilde\beta}(V,T)}\Big(
			(C_1+C_2 \tau^\pm_{\widetilde\beta}(V,T)^{p})^q+
            (C_1+C_2( \tau^\pm_{\widetilde\beta}(V,T)+C)^{p})^q\Big),
        \end{equation*}
        and
        $$\overline \Lambda_T = \sup_{T'\geq T}\sup_{0\leq t < T'} \Lambda(T';V,t).$$
	\begin{proof}
        
        Note that if $t\leq \theta_{\widetilde\beta}^-(V,T)$ then $ \tau_\beta(V,t) \leq \theta_{\widetilde\beta}^-(V,T) - t \leq T-t$, which implies $\tau_\beta(V,t)=\widehat\tau_\beta(T;V,t)$. We also get, for any $ t\in\N $ with $ \theta_{\widetilde\beta}^-(V,T) < t \leq \theta_{\widetilde\beta}^+(V,T) $, that $\tau_\beta(V,t) \leq \theta_{\widetilde\beta}^+(V,T) - t \leq  \tau^\pm_{\widetilde\beta}(V,T)$.

        This yields
        \begin{align*}
            &\sum_{t=0}^{T-1}\big|f(\tau_\beta(V,t)) - f(\widehat{\tau}_\beta(T;V,t))\big|^q\\
			&=\sum_{t=0}^{T-1}\1_{\tau_\beta(t) \geq T-t}
			\big|f(\tau_\beta(V,t))-f(T-t+\Lambda(T;V,t))\big|^q\\
			&\leq\sum_{t=(\theta_{\widetilde\beta}^-(V,T)+1) \lor 0}^{\theta_{\widetilde\beta}^+(V,T)\land (T-1)}
			\big|f(\tau_\beta(V,t))-f(T-t+\Lambda(T;V,t))\big|^q\\
			&\leq\sum_{t=(\theta_{\widetilde\beta}^-(V,T)+1) \lor 0}^{\theta_{\widetilde\beta}^+(V,T)\land (T-1)}
			\left(\big|C_1+C_2 \tau^\pm_{\widetilde\beta}(V,T)^{p}\big|^q+\big|C_1+C_2( \tau^\pm_{\widetilde\beta}(V,T)+\Lambda(T;V,t))^{p}\big|^q\right)\\
			&\leq \tau^\pm_{\widetilde\beta}(V,T)\Big(
			(C_1+C_2 \tau^\pm_{\widetilde\beta}(V,T)^{p})^q+
            (C_1+C_2( \tau^\pm_{\widetilde\beta}(V,T)+\overline \Lambda_T)^{p})^q\Big).
		\end{align*}		
	\end{proof}
\end{lemm}

Since $\Pp_{p,q}(T,C)$ is asymptotically a power function of $\tau^\pm_{\widetilde\beta}(V,T)$ we can guarantee different forms of convergence under moment conditions on $\tau^\pm_{\widetilde\beta}(V,0)$. To simplify our assumptions we first prove that finite moments $\tau^\pm_{\widetilde\beta}(V,0)$ is equivalent to moments of $\tau_{\widetilde\beta}(V,0)$.

\begin{lemm}\label{lem:transitmoment}
    Assume $\widetilde\beta\subset\R^N$. 

    If, for any $\alpha\in[0,\infty)$, we have $\E[\tau_{\widetilde{\beta}}(V,t)^\alpha]<\infty$, then $\E[\tau^\pm_{\widetilde\beta}(V,t)^\alpha]<\infty$.
    \begin{proof}
        First denote the backwards hitting time $\tau_{\widetilde\beta}^{-}(V,t) \triangleq \min\{s\geq 0:V_{t-s} \in \beta_s\}$.
        
        Further define $H^i_j=\{\omega: \tau_{\widetilde\beta}^{-}(V,t)=i, \, \tau_{\widetilde\beta}(V,t)=j\}$, and note that by ergodicity we have $\PP(H^i_j)=\PP(S^{i-1}(H^i_j))=\PP(H^1_{j+i-1})$ for $i,j\geq 1$. The disjoint sets $\{H^i_j\}_{i,j\in\N}$ forms the basis of the \textit{Kakutani skyscraper}, and are commonly applied in relation to transit times of ergodic processes (see e.g. \cite{ergodicbook}).
        
        To simplify notation, we will introduce $k=i+j$ as an index variable. For any $\alpha\in[0,\infty)$, this yields
        \begin{align*}
            \E[\tau_{\widetilde\beta}(V,t)^\alpha]
            &= \sum_{i=1}^\infty \sum_{j=1}^\infty \PP(H^i_j) j^\alpha
            = \sum_{k=2}^\infty \sum_{j=1}^{k-1} \PP(H^i_j) j^\alpha\\
            &= \sum_{k=2}^\infty \sum_{j=1}^{k-1} \PP(H^1_{k-1}) j^\alpha
            = \sum_{k=2}^\infty \PP(H^1_{k-1}) \sum_{j=1}^{k-1}  j^\alpha\\
            &\geq \sum_{k=2}^\infty  \PP(H^1_{k-1}) \int_0^{k-1} j^\alpha dj
            = \sum_{k=2}^\infty  \PP(H^1_{k-1}) \frac{(k-1)^{\alpha+1}}{\alpha+1}.\\
        \end{align*}
        Similarly, we have
        \begin{align*}
            \E[\tau^\pm_{\widetilde\beta}(V,t)^\alpha]
            &= \sum_{i=1}^\infty \sum_{j=1}^{\infty} \PP(H^i_j) {(i+j-1)}^\alpha
            = \sum_{k=2}^\infty \sum_{j=1}^{k-1} \PP(H^i_j) (i+j-1)^\alpha\\
            &= \sum_{k=2}^\infty \sum_{j=1}^{k-1} \PP(H^1_{k-1}) (k-1)^\alpha
            = \sum_{k=2}^\infty \PP(H^1_{k-1}) \sum_{j=1}^{k-1}  (k-1)^\alpha\\
            &= \sum_{k=2}^\infty  \PP(H^1_{k-1}) (k-1)^{\alpha+1}
            \leq  (\alpha+1)\E[\tau_{\widetilde\beta}(V,t)^\alpha]<\infty.
        \end{align*}
    \end{proof}
\end{lemm}

With this we can control the error arising from replacing $\tau_{\widetilde\beta}(V,t)$ by $\widehat\tau_{\widetilde\beta}(V,t)$.

\begin{lemm}\label{exceedfeasiblefull}
	%
	%
	Assume the threshold function, $ \beta $, satisfies $ \beta_s\supseteq\widetilde{\beta} $ for some measurable $\widetilde\beta\subseteq\R^N$ and all $ s\geq 0 $ with $ \E\left[\tau_{\widetilde{\beta}}^{pq+1} \right] <\infty$ for some $ p,q \in [0,\infty) $. 
	
	Fix the functions $ f:\N\mapsto\R $ and $\Lambda: \N\times(\R^N)^\Z\times\N\mapsto\R $ and 
    assume that
	$$ \limsup_{T}\sup_{0\leq t  < T}|\Lambda(T;V,t)| < \infty, \quad |f(t)| \leq C_1 + C_2t^{p}, $$
	almost surely for some $ C_1,\,C_2\in[0,\infty) $. 
 
	Under these assumptions we have that
	$$ \lim_{T\to\infty} \frac{1}{T}\sum_{t=0}^{T-1}\big|f(\tau_\beta(V,t)) - f(\widehat{\tau}_\beta(T;V,t))\big|^q=0,$$
	almost surely. 
	\begin{proof}

        We start by setting $ C_3 $ to be any number such that $ C_3 > \limsup_{T}|\Lambda(T;V,t)|$ and $ C_3\in \Q $. We may then assume without loss of generality that $T$ is sufficiently large to ensure  $C_3 \geq \sup_{T'\geq T}\sup_{0\leq t<T'}|\Lambda(T';V,t)|$. By \Cref{lem:bound_by_P}, we have 
        
        $$ \frac{1}{T}\sum_{t=0}^{T-1}\big|f(\tau_\beta(V,t)) - f(\widehat{\tau}_\beta(T;V,t))\big|^q
        \leq \frac{1}{T}\Pp_{p,q}(T,C_3).$$
        We note, for any deterministic $ C $, that $\Pp_{p,q}(T,C)$ has a finite expectation due to our assuption of $ \E\left[\tau_{\widetilde{\beta}}^{pq+1} \right] <\infty$ and \Cref{lem:transitmoment}. This implies, by \Cref{thm:mainergodic}, that
		\begin{align*}
			&\lim_{T\to\infty}\frac{1}{T}\Pp_{p,q}(T,C)\\
			&= \lim_{T\to\infty} \left(\frac{T+1}{T}\right)\frac{1}{T+1}
			\sum_{t=0}^{T}\Pp_{p,q}(t,C)
			-\lim_{T\to\infty} \frac{1}{T} \sum_{t=0}^{T-1}\Pp_{p,q}(t,C)\\
			&=\E\left[\Pp_{p,q}(0,C)\right]-\E\left[\Pp_{p,q}(0,C)\right]\\
			&=0,
		\end{align*}
		
		\noindent almost surely. This further implies that
		
		\begin{align*}
			&\PP\left(\lim_{T\to\infty}\frac{1}{T}\Pp_{p,q}(T,C)=0 \text{ for all } C\in\Q \right)\\
			&=\PP\left(\bigcap_{C\in\Q}
			\left\{
			\lim_{T\to\infty}\frac{1}{T}\Pp_{p,q}(T,C)=0
			\right\} \right)\\
            &=1-\PP\left(\bigcup_{C\in\Q}
			\left\{
			\lim_{T\to\infty}\frac{1}{T}\Pp_{p,q}(T,C) \neq 0
			\right\} \right)\\
			&\geq1-\sum_{C\in\Q} \PP\left(
			\left\{
			\lim_{T\to\infty}\frac{1}{T}\Pp_{p,q}(T,C) \neq 0
			\right\} \right)\\
			&=1.
		\end{align*}

		\noindent This ensures that $ \Pp_{p,q}(T,C_3)/T $ converges almost surely to $ 0 $, and therefore that $\lim_{T\to\infty} \sum_{t=0}^{T-1}\big|f(\tau_\beta(V,t)) - f(\widehat{\tau}_\beta(T;V,t))\big|^q/T=0$ almost surely.
		
	\end{proof}
\end{lemm}

\begin{rem}\label{rem:nomoment}
    In the case where $ f(\tau)=\1(\tau> n) $ we can omit the moment condition on $ \tau_{\widetilde{\beta}} $, as well as $ \beta_s\supseteq \widetilde{\beta}$ and $\limsup_{T}\sup_{0\leq t  < T}|\Lambda(T;V,t)| < \infty$. This is done by bounding $\sum_{t=0}^{T-1}\big|f(\tau_\beta(V,t)) - f(\widehat{\tau}_\beta(T;V,t))\big|^q$ by $n2^q/T$ in the proof of \Cref{lem:bound_by_P}.
\end{rem}

With \Cref{exceedfeasiblefull} we can guarantee that our estimator converges almost surely to the desired moment. This result also allows for the computation of the distribution of $f(\tau_\beta(V,0))$, along with its correlation structure.

\begin{thm}\label{thm:RP_est}
    Assume that the threshold function $ \beta $ satisfies $ \beta_s\supseteq\widetilde{\beta} $ for some measurable $\widetilde\beta\subseteq\R^N$ and all $ s\geq 0 $. 
	
	Fix the functions $ f:\N\mapsto\R $ and $\Lambda: \N\times(\R^N)^\Z\times\N\mapsto\R $ and 
    assume that
	$$ \limsup_{T}\sup_{0\leq t  < T}|\Lambda(T;V,t)| < \infty, \quad |f(\tau)| \leq C_1 + C_2\tau^{p}, $$
    almost surely for some $ C_1,\,C_2\in[0,\infty) $. Under these conditions we have the following.

    If $\E\left[\tau_{\widetilde{\beta}} \right] <\infty$, then 
    $$\PP\left(f(\tau_\beta(V,0))> s\right)  
		= \lim_{T\to\infty}\frac{1}{T}\sum_{t=0}^{T-1}\1(f(\widehat\tau_\beta(T;v,t)> s)).$$

    If $\E\left[\tau_{\widetilde{\beta}}^{p+1} \right] <\infty$, then 
    $$\E\left[f(\tau_\beta(V,0)) \right]  
		= \lim_{T\to\infty}\frac{1}{T}\sum_{t=0}^{T-1}f(\widehat\tau_\beta(T;V,t)).$$

    If $\E\left[\tau_{\widetilde{\beta}}^{2p+1} \right] <\infty$, then 
    $$\E\left[f(\tau_\beta(V,0))f(\tau_\beta(V,i)) \right]  
		= \lim_{T\to\infty}\frac{1}{T}\sum_{t=i}^{T-1}f(\widehat \tau_\beta(T;V,t-i))f(\widehat\tau_\beta(T;V,t)).$$
    \begin{proof}
        The first two equations are immediate consequences of \Cref{exceedfeasiblefull} with $q=1$. For the third equation, we define $a_{t,T}=f(\widehat\tau_\beta(T;V,t))$ and $\alpha_t=f(\tau_\beta(V,t))$ as shorthand notations. We then get by the Cauchy–Schwarz inequality that
        \begin{align*}
            &\left|\E\left[f(\tau_\beta(V,0))f(\tau_\beta(V,i)) \right]  
		- \lim_{T\to\infty}\frac{1}{T}\sum_{t=i}^{T-1}a_{t,T}a_{t-i,T} \right|\\  
		&\leq \lim_{T\to\infty}\frac{1}{T}\sum_{t=i}^{T-1}\left|a_{t,T}a_{t-i,T}-\alpha_t\alpha_{t-i}  \right|\\
        &\leq \lim_{T\to\infty}\frac{1}{T}\sum_{t=i}^{T-1}\Big(\left|a_{t,T}(a_{t-i,T}-\alpha_{t-i})\right|+\left|(a_{t,T}-\alpha_t)\alpha_{t-1}\right|  \Big)\\
        &\leq \lim_{T\to\infty}\sqrt{\frac{1}{T}\sum_{t=i}^{T-1}a_{t,T}^2}\sqrt{\frac{1}{T}\sum_{t=i}^{T-1}(a_{t-i,T}-\alpha_{t-i})^2}\\
        & \quad +\lim_{T\to\infty}\sqrt{\frac{1}{T}\sum_{t=i}^{T-1}\alpha_{t-i}^2}\sqrt{\frac{1}{T}\sum_{t=i}^{T-1}(a_{t,T}-\alpha_t)^2}.        
        \end{align*}
        By \Cref{thm:mainergodic} and \Cref{exceedfeasiblefull} (for $q=2$), this limit is 0.
        
    \end{proof}
\end{thm}

\begin{rem}\label{rem:nomomentfinal}
    As per \Cref{rem:nomoment}, we can simplify the requirements for the first equality if $f(\tau)=\tau$. We may omit the conditions $\E\left[\tau_{\widetilde{\beta}} \right] <\infty$,  $ \beta_s\supseteq \widetilde{\beta}$, and $\limsup_{T}\sup_{0\leq t  < T}|\Lambda(T;V,t)| < \infty$.
\end{rem}

Now that we have proven our law of large numbers, we can move on to our last major result, the central limit theorem for our adjusted estimator.

\begin{thm}\label{thm:feasibleCLT}
    Assume the threshold function, $ \beta $, satisfies $ \beta_s\supseteq\widetilde{\beta} $ for some measurable $\widetilde\beta\subseteq\R^N$ and all $ s\geq 0 $ with $ \E\left[\tau_{\widetilde{\beta}}^{\max(2p,1+p)} \right] <\infty$ for some $ p\in [0,\infty) $. 
	
	Fix the functions $ f:\N\mapsto\R $ and $\Lambda: \N\times(\R^N)^\Z\times\N\mapsto\R $ and 
    assume that
	$$ \limsup_{T}\sup_{0\leq t  < T}|\Lambda(T;V,t)| < \infty, \quad |f(t)| \leq C_1 + C_2t^{p}, $$
	almost surely for some $ C_1,\,C_2\in[0,\infty) $. 
    
    Define $\M_T=\sum_{t=0}^{T-1}f(\tau_\beta(V,t))$, $\widehat \M_T=\sum_{t=0}^{T-1}f(\widehat\tau_\beta(T;V,t))$, and $\sigma_T^2=\V ar[\M_T^2]$. Note that $ \E\left[\tau_{\widetilde{\beta}}^{2p} \right] <\infty$ implies that $\sigma_T<\infty$.
    
    Further assume that $t \mapsto \tau_\beta(V,t)$ satisfies the remaining conditions of \Cref{thm:mainCLT}. This entails that $t \mapsto \tau_\beta(V,t)$ is a strongly mixing process with $\lim_{T\to\infty}\sigma_T=\infty$, and that the family $\{\M_T^2/\sigma_T^2\}_{T\in\N}$ is uniformly integrable. 
    
    We then have that $(\widehat \M_T-\E[\M_T])/\sigma_T$ converges in distribution to a standard normal random variable.

    \begin{proof}    
        From \Cref{thm:mainCLT} we know that $( \M_T-\E[\M_T])/\sigma_T$ converges in distribution to a standard normal random variable. 
        
        We also have the following standard convergence result for sequences of random variables (see e.g.\ \cite{van2000asymptotic}). If two sequences, $\{A_n\}_{n\in\N},\, \{B_n\}_{n\in\N} $ satisfies $A_n\to A$ in distribution and $|A_n-B_n|\to 0 $ in probability, then $B_n\to A$ in distribution.

        Therefore it suffices to show that $|\M_T-\widehat \M_T|/\sigma_T \to 0$ in probability. By \Cref{lem:bound_by_P} we have
        $$\frac{|\M_T-\widehat \M_T|}{\sigma_T}\leq \frac{\Pp_{p,1}(T,\overline \Lambda_T)}{\sigma_T},$$
        where $\overline \Lambda_T = \sup_{T'\geq T}\sup_{0\leq t < T'} \Lambda(T';V,t)$.

        Since $\limsup_{T\to\infty}\overline \Lambda_T<\infty$ almost surely, we have $\lim_{T\to\infty}\PP(\overline \Lambda_T=\infty)=0$. Consequently, for any $\delta>0$, there must be some $T_1\in\N$ such that $\PP( \overline \Lambda_{T_1} = \infty) < \delta/3$. We can then further pick a $C_\text{min}\in[0,\infty)$ such that $\PP( \overline \Lambda_{T_1} > C_\text{min}) < \delta/2$. Next, by the fact that $T\mapsto \overline \Lambda_T$ is monotone non-increasing, we have
        $$\PP( \overline \Lambda_T > C_\text{min}) < \delta/2,$$
        for all $T\geq T_1$. 

        We further note, for any deterministic $C$, that $T\mapsto \Pp_{p,1}(T,C)$ is a finite and stationary sequence. Consequently, since $\sigma_T\to\infty$, we get ${\Pp_{p,1}(T,C)}/{\sigma_T} \to 0 $ in probability. As such, there must be some $T_2\in\N$ such that 
        $$\PP\left(\frac{\Pp_{p,1}(T,C_\text{min})}{\sigma_T}>\epsilon\right)<\delta/2,$$
        for all $T\geq T_2$.

        Finally, we note that $C \mapsto \Pp_{p,1}(T,C)$ is monotone non-decreasing, which implies 
        $$\Pp_{p,1}(T,C_\text{min}) \geq \Pp_{p,1}(T,\overline \Lambda_T),$$
        conditional on $\overline \Lambda_T \leq C_\text{min}$.
        
        Combining these equations we get
        \begin{align*}
            \PP\left({\Pp_{p,1}(T,\overline \Lambda_T)}/{\sigma_T}>\epsilon\right)
            =\,&
            \PP\left(\Pp_{p,1}(T,\overline \Lambda_T)/\sigma_T>\epsilon \big| \overline \Lambda_T > C_\text{min} \right)   
                \PP\left( \overline \Lambda_T > C_\text{min}\right)\\
            &+\PP\left(\Pp_{p,1}(T,\overline \Lambda_T)/\sigma_T>\epsilon \big| \overline \Lambda_T \leq C_\text{min}\right)   
                \PP\left( \overline \Lambda_T \leq C_\text{min}\right)\\
            \leq\,& 
            \PP\left( \overline \Lambda_T > C_\text{min}\right)\\
            &+\PP\left(\{\Pp_{p,1}(T,\overline \Lambda_T)/\sigma_T>\epsilon\} \bigcap \{\overline \Lambda_T \leq C_\text{min}\}\right)\\
            \leq\,& 
            \delta/2
            +\PP\left(\{\Pp_{p,1}(T,C_\text{min})/\sigma_T>\epsilon\} \bigcap \{\overline \Lambda_T \leq C_\text{min}\}\right)\\
            \leq\,& 
            \delta/2
            +\PP\left(\Pp_{p,1}(T,C_\text{min})/\sigma_T>\epsilon \right)\\
            <\,& \delta,
        \end{align*}
        for all $T\geq \max(T_1,T_2)$. As such we have ${\Pp_{p,1}(T,\overline \Lambda_T)}/{\sigma_T} \to 0 $, in probability, which implies $(\widehat \M_T - T\E[f(\tau(V,0))])/\sigma_T$ converges, in distribution, to a standard normal random variable.
        
    \end{proof}
\end{thm}

Lastly we prove the general form of \eqref{eq:statmean}, i.e.\ our simplified formula for our estimator in the case where $\beta$ is constant.

\begin{prop}\label{prop:stat_mom}
	Let $ \beta_s = \widetilde{\beta} $ for all $ s $ and define, $ \Theta $, the set of all observed points where the process is outside $ \widetilde{\beta} $, i.e. $ \Theta\triangleq\{\theta\in\N: \theta<T,\, V_\theta \in\widetilde{\beta}\} $. We assume that the size of this set, $ |\Theta| $, is at least $ 1 $. We denote by $ \theta_i $, $ i\leq|\Theta| $ the $ i $'th element of $ \Theta $, and define $ \theta_{|\Theta|+1}=T+\theta_1 $. We lastly introduce the shorthand $\Delta\theta_i=\theta_{i+1}-\theta_i$. This yields, for any $f:\N\to \R$, that
	
	$$\frac{1}{T}\sum_{t=0}^{T-1} \widehat\tau_\beta(t)^k
	=\sum_{i=1}^{|\Theta|}
            \sum_{j=0}^{\Delta\theta_i-1}
            \frac{f(j)}{T},$$
	where we again use
	$$\widehat\tau_\beta(t)=\min\{s\geq 0:V_{t+s \text{ mod } T}\in\beta_s\},
	\quad a \text{ mod } b = a-b\lfloor a/b\rfloor.$$
	
	\begin{proof}
		We note that for all $ t $ we have $ t\leq\theta_{i(t)} $ for some minimal $ i(t) $, which implies $ \widehat\tau_\beta(t) =  \theta_{i(t)}-t$. This means that
		
		\begin{align*}
			&\frac{1}{T}\sum_{t=0}^{T-1} f(\widehat\tau_\beta(t))\\
            =&\frac{1}{T}\left(
			\sum_{t=0}^{\theta_{1}} f(\widehat\tau_\beta(t))
			+\sum_{i=1}^{|\Theta|-1}\sum_{t=\theta_{i}+1}^{\theta_{i+1}} f(\widehat\tau_\beta(t))
			+\sum_{t=\theta_{|\Theta|}}^{T-1} f(\widehat\tau_\beta(t))
			\right)\\
			=&\frac{1}{T}\left(
			\sum_{t=0}^{\theta_{1}} f\left(\theta_{1}-t\right)
			+\sum_{i=1}^{|\Theta|-1}\sum_{t=\theta_{i}+1}^{\theta_{i+1}} f\left(\theta_{i+1}-t\right)
			+\sum_{t=\theta_{|\Theta|}}^{T-1} f\left(\theta_{1}+T-t\right)
			\right)\\
			=&\frac{1}{T}\left(
			\sum_{t=T}^{\theta_{1}+T} f\left(\theta_{1}+T-t\right)
			+\sum_{i=1}^{|\Theta|-1}\sum_{t=\theta_{i}+1}^{\theta_{i+1}} f\left(\theta_{i+1}-t\right)
			+\sum_{t=\theta_{|\Theta|}}^{T-1} f\left(\theta_{1}+T-t\right)
			\right)\\
			=&\frac{1}{T}\left(
			\sum_{i=1}^{|\Theta|}\sum_{t=\theta_{i}+1}^{\theta_{i+1}} f\left(\theta_{i+1}-t\right)
			\right)\\
			=&
			\sum_{i=1}^{|\Theta|}
            \sum_{j=0}^{\Delta\theta_i-1}
            \frac{f(j)}{T}.
		\end{align*}
	\end{proof}
\end{prop}

The most common moments of interest are $\tau_\beta(V,0)^p$ for $p=1$ and $p=2$. As such we also give a more explicit form for these cases.
\begin{cor}
    Let the conditions of \Cref{prop:stat_mom} hold. We the have
    $$\frac{1}{T}\sum_{t=0}^{T-1} \widehat\tau_\beta(t)
	=\sum_{i=1}^{|\Theta|}
            \frac{\Delta\theta_i(\Delta\theta_i-1)}{2 T},$$
    $$\frac{1}{T}\sum_{t=0}^{T-1} \widehat\tau_\beta(t)^2
	=\sum_{i=1}^{|\Theta|}
            \frac{\Delta\theta_i(\Delta\theta_i-1)(2\Delta\theta_i-1)}{6 T}.$$
\end{cor}

\section*{Acknowledgements}
The author acknowledge financial support by the Research Council of Norway under the SCROLLER project, project number 299897.

\bibliographystyle{abbrvdin}
\bibliography{main.bib}

\end{document}